\documentclass{amsart}
\usepackage{rlepsf, graphicx, amssymb, epstopdf, amsmath, amssymb}

\theoremstyle{plain}
\newtheorem{Thm}{Theorem}
\newtheorem{Lem}[Thm]{Lemma}
\newtheorem{Cor}[Thm]{Corollary}
\newtheorem{Prop}[Thm]{Proposition}

\theoremstyle{definition}

\theoremstyle{remark}

\def\C{\mathbb{C}}
\def\R{\mathbb{R}}
\def\Z{\mathbb{Z}}

\def\sx{\hbox{$S^2 \times S^2$}}
\def\snx{\hbox{$S^2 \widetilde{\times} S^2$}}

\def\bdy{\partial}

\def\Lf{Lefschetz fibration }
\def\aLf{achiral Lefschetz fibration }

\def\dfn#1{{\em #1}}

\title{Realizing $4$-manifolds as achiral Lefschetz fibrations}
\author{John B.\ Etnyre}
\address{
    School of Mathematics,
    Georgia Institute of Technology,
    686 Cherry St.,
    Atlanta, GA  30332-0160}
\email{etnyre@math.gatech.edu}
\urladdr{http://math.gatech.edu/\char126 etnyre}

\author{Terry Fuller}
\address{Department of Mathematics,
California State University, Northridge, 
Northridge, CA 91330 } 
\email{terry.fuller@csun.edu}
\urladdr{http://www.csun.edu/\char126 tf54692}

\begin{document}
\begin{abstract}
We show that any 4-manifold, after surgery on a curve, admits an achiral Lefschetz fibration. In particular, if $X$ is a
simply connected 4-manifold we show that $X\# S^2\times S^2$ and $X\# S^2\widetilde\times S^2$ both admit achiral Lefschetz 
fibrations. We also show these surgered manifolds admit near-symplectic structures and prove more generally that achiral Lefschetz fibrations with sections
have near-symplectic structures. 
As a corollary to our proof we obtain an alternate proof of Gompf's result on the existence of 
symplectic structures on Lefschetz pencils.
\end{abstract}
\maketitle

\section{Introduction}
Symplectic $4$-dimensional manifolds are known to be characterized as those admitting the structure of a Lefschetz fibration. More precisely, Donaldson \cite{Donaldson99} proved that every symplectic $4$-manifold admits a Lefschetz pencil, which can be blown up at its base points to yield a Lefschetz fibration. Conversely, Gompf \cite{GompfStipsicz99} showed that any $4$-manifold with a \Lf admits a symplectic structure, provided the fibers are non-trivial in homology.

The definition of a \Lf includes the provision that the orientations in the local holomorphic description of a critical point match the global orientations of the total and base spaces of the fibration. This condition is crucial for the above results, for symplectic structures serve to orient the manifolds involved, and Donaldson and Gompf each elucidate how symplectic structures on fibers are compatible with a global symplectic structure. Therefore if one relaxes this requirement, the resulting wider class of fibrations, known as {\it achiral} Lefschetz fibrations, will no longer respect symplectic structures. It is natural to ask which arbitrary (i.e. not necessarily symplectic) smooth manifolds admit achiral Lefschetz fibrations.

The first result concerning the existence of achiral Lefschetz fibrations is due to Harer \cite{Harer79} who proved that a
4-manifold that has a handle decomposition with only a 0-handle, 1-handles and 2-handles admits an achiral 
Lefschetz fibration over the disk. It was observed in \cite{GompfStipsicz99} that any closed simply-connected 4-manifold admits 
an achiral Lefschetz fibration over $S^2$ after connect summing with $S^2\times S^2$ some number of
times--this number is unknown and depends on the manifold. There is no similar statement for non-simply
connected 4-manifolds. 

In the other direction, the only known obstruction to the existence of an achiral Lefschetz 
fibration, also found in \cite{GompfStipsicz99}, is that for
a manifold with positive definite intersection form the inequality
\[1-b_1+b_2\geq q\geq 0.\]
must hold, where $q$ is the number of negative vanishing cycles. (There is an analogous result for
negative definite manifolds.) Thus, for example, $\#_n S^1\times S^3$ does not admit an achiral Lefschetz fibration
for $n>1.$

Our main result is the following.

\begin{Thm} \label{main} Let $X$ be a smooth, closed, oriented $4$-manifold. Then there exists a framed circle in
$X$ such that the manifold obtained by surgery along that circle admits an achiral Lefschetz
fibration with base $S^2$. Moreover, all these fibrations admit sections.
\end{Thm}

Since surgery on a circle in simply-connected $4$-manifolds always changes the manifold by a
connected sum with an $S^2$-bundle over $S^2$, we immediately see that $X\# \sx$ or $X\#
\snx$ admits an achiral Lefschetz fibration whenever $X$ is simply-connected. This
can be strengthened to the following result.

\begin{Cor}\label{cor1} Let $X$ be a smooth, closed, simply-connected $4$-manifold. Then both $X\# \sx$ and $X\#
\snx$ admit an achiral Lefschetz fibration.
\end{Cor}

Recently work of Taubes \cite{Taubes98} has created a great deal of interest in near-symplectic structures. These are
closed 2-forms on a 4-manifold that are symplectic off of an embedded 1-manifold, and vanish in a prescribed way
along this 1-manifold (see Section~\ref{nss}). We prove the following about achiral Lefschetz structures and near-symplectic
structures.
\begin{Thm}\label{nssonlf}
If a 4-manifold admits an achiral Lefschetz fibration over $S^2$ with a section, then it has an near-symplectic structure.
Moreover, the near-symplectic structure can be chosen so that any pre-assigned disjoint sections are symplectic.
\end{Thm}
We note that since $b_2^+>0$ for the achiral Lefschetz fibrations of Theorem~\ref{main}, the existence of a near-symplectic structure follows from a result of Honda \cite{Honda04a}.
Our proof, however, gives a more explicit construction of the near-symplectic form ({\em cf.} \cite{GayKirby04}) and 
illuminates the relationship between the achiral Lefschetz structure and the near-symplectic structure. This allows, among other things, for the
possibility of a Donaldson-Smith approach to studying symplectic submanifolds/holomorphic curves in these manifolds, as in \cite{DonaldsonSmith04}.
Combining this theorem with Theorem~\ref{main} yields the following result.
\begin{Cor}
Let $X$ be a smooth, closed, oriented $4$-manifold. Then there exists a framed circle in
$X$ such that the manifold obtained by surgery along that circle admits a near-symplectic structure.
Moreover, if $X$ is simply connected then both $X\# S^2\times S^2$ and $X\# S^2\widetilde{\times} S^2$ admit
a near-symplectic structure. 
\end{Cor}
The method of proof for Theorem~\ref{nssonlf} yields a different proof of the well-known result of Gompf mentioned above.
\begin{Thm}[Gompf, \cite{GompfStipsicz99}]
If a 4-manifold $X$ admits a Lefschetz fibration over $S^2$ with a section, then it has a symplectic structure. Moreover,
the symplectic structure may be chosen so that any preassigned disjoint sections are symplectic.
\end{Thm}
This 
is actually weaker that Gompf's result, where one does not need to assume the existence of a section, 
only that the fiber is non-trivial in homology; however, we are still able to recover an important corollary of Gompf's result. 
\begin{Cor}
If a 4-manifold $X$ admits a Lefschetz pencil, then it admits a symplectic structure.
\end{Cor}
This corollary was previously observed to follow from arguments similar to ours by David Gay \cite{Gay02}.

{\it Acknowledgments}: The first author was partially supported by NSF Career Grant (DMS-0239600) and FRG-0244663.

\section{Lefschetz fibrations, open book decompositions and handlebodies}\label{sec:lf}
A \dfn{Lefschetz fibration} of an oriented 4-manifold $X$ is a map $f\colon\, X\to F$ to a surface $F$ such that 
all the critical points of $f$ lie in the interior of $X$ and for each critical point there is
an orientation preserving  coordinate chart on which  $f\colon \C^2\to \C$ takes the 
form $f(z_1,z_2)=z_1z_2.$ We assume all the critical points occur on distinct fibers. 

If $x$ is a non-critical value in $F$ then $\Sigma=f^{-1}(x)$ is a surface properly embedded in $X.$ The diffeomorphism type of $f^{-1}(x)$ is
independent of the non-critical value $x$, and may have boundary, if $X$ does. Let $p$ be a critical point in $X$ and $U$ a closed 
disk neighborhood of $f(p)$
in $F$ that contains no other critical values. If $y\in\partial U$ and $c$ is a radial path in $U$ from $y$ to $f(p)$ 
then there is an embedded disk $D_c$ in $X$ that projects to $\gamma$ and $f^{-1}(x)\cap D_c$
is a simple closed curve $\gamma_p$ in the fiber above $x$ for all $x\in(c \setminus \{f(p)\}).$ 
(Note we use $\gamma_p$ for all curves in the fibers above $c$.) Note that $\gamma_p$ will usually be 
non-trivial in the homology of the fiber but will be trivial in the homology of $f^{-1}(U).$
The curve $\gamma_p$ is called the \dfn{vanishing cycle} associated to $p.$
It can be shown that $f^{-1}(U)$ is obtained from $\Sigma\times D^2$ by attaching a 2-handle to $\gamma_p$ with framing
one less than the framing induced on $\gamma_p$ by $f^{-1}(y).$ In addition $f^{-1}(\partial U)$ is a $\Sigma$-bundle 
over $S^1=\partial U$ with monodromy given by a positive Dehn twist along $\gamma_p,$ which we denote 
$D_{\gamma_p}.$

More generally, if $F=D^2$ then we fix a point $y\in\partial D^2$ and a collection of embedded arcs $c_1,\ldots, c_k$ connecting
$y$ to the critical points $p_1,\ldots, p_k,$ such that they only intersect at $y.$ We order the $c_i$'s so that a small circle about
$y$ intersects them in a counter-clockwise order. We now have a collection of vanishing
cycles $\gamma_{p_1},\ldots, \gamma_{p_k}$ in $\Sigma=f^{-1}(y).$ The manifold $X$ is obtained from $\Sigma\times D^2$ by
attaching 2-handles along the $\gamma_{p_i}$'s with framing one less than the framing induced by $\Sigma.$ Moreover, 
$f^{-1}(\partial D^2)$ is a $\Sigma$-bundle over $S^1$ with monodromy $D_{p_1}\circ\dots\circ D_{p_k}.$

If the fibers of the Lefschetz fibration do not have boundary (and $F=D^2$) then $\partial X=f^{-1}(\partial D^2)$
is the surface bundle described above. If the fibers do have boundary then $\partial X=((\partial \Sigma)\times D^2)
\coprod (f^{-1}(\partial D^2)).$ Clearly $(\partial \Sigma)\times D^2$ is a neighborhood of $B=\partial f^{-1}(x)\subset 
\partial X$ for any $x$ in the interior of $D^2.$ Thus it is easy to see that $(\partial X)\setminus B$ is the
$\Sigma$-bundle over $S^1$ described above. More specifically, the Lefschetz fibrations induces an \dfn{open book
decomposition} of $\partial X.$ Recall an \dfn{open book decomposition} of a 3-manifold $M$ is a pair $(B,\pi)$ where
$B$ is an oriented link in $M$ and $\pi\colon\, M\setminus B \to S^1$ is a fibration of the complement of $B$ such 
that $\partial \overline{\pi^{-1}(\theta)}=B$ for all $\theta\in S^1.$ The fibers of $\pi$ are called \dfn{pages} of
the open book and $B$ is called the \dfn{binding} of the open book. For a more complete discussion of topological 
Lefschetz fibrations and open book decompositions see \cite{EtnyreOBN, GompfStipsicz99}. 

An \dfn{achiral Lefschetz fibration} is an oriented 4-manifold $X$ and a map $f\colon\, X\to F$ exactly as in
the definition of Lefschetz fibration above except that the coordinate charts do not have to be orientation preserving.
Critical points with non-orientation preserving charts will be called a \dfn{negative} critical point. 
The entire discussion above carries over to the achiral case, except that the 2-handles attached to the
vanishing cycle of a negative critical point will have framing one more than the framing induced by $\Sigma$ and
the contribution to the monodromy will be a left handed Dehn twist $D^{-1}_{\gamma_p}$ about the vanishing cycle. 

A key theorem for the proof of our main theorem is the following result.
\begin{Thm}[Harer 1979, \cite{Harer79}]\label{thm:harer}
Let $X$ be a 4-dimensional handlebody with all handles of index less than or equal to two. Then $X$ admits an achiral Lefschetz fibration over $D^2$ with bounded fibers. 
\end{Thm}

\begin{proof}[Sketch of Proof]
We briefly sketch Harer's original proof, since it is very nice and not easy to find in the literature
({\em cf.} \cite{AkbulutOzbagci01}). Recall if $\Sigma$ is any oriented surface with boundary, then $\Sigma\times D^2$ is diffeomorphic to 
$B^4$ with $k$ 1-handles attached, where $k=-\chi(\Sigma)+1.$  Thus if $X$ has $k$ 1-handles,  then we begin by letting $\Sigma$ be a disk with $k$ open 
disks removed. We can picture $\Sigma$ as a disk-with-bands inside of the 1-handlebody $X_1$ of $X$ as on the left side of Figure~\ref{harer1}.
\begin{figure}[ht]
 \relabelbox \small {
  \centerline{\epsfbox{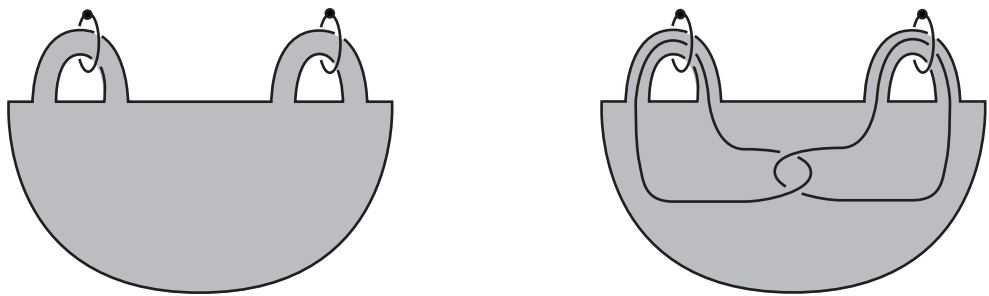}}}
  \endrelabelbox
        \caption{The surface $\Sigma$ in $X_1.$}
        \label{harer1}
\end{figure}
Note that we have constructed a trivial Lefschetz structure for  $X_1.$ 

We next consider the attaching link $L$ for the 2-handles of $X$. We first isotope $L$ into a neighborhood of $\Sigma$
so that it projects onto $\Sigma$ with only double points, as on the right side of Figure~\ref{harer1}. To prove the theorem we must modify $\Sigma$ so that the attaching circle of each 2-handle can be embedded  on a separate fiber of $\partial(\Sigma\times D^2)$, and arrange that each is attached with framing $\pm 1$ relative to the product framing on $\Sigma\times D^2$. Ignoring the framings momentarily, we can accomplish the former by forming the connected sum of $\Sigma$ with a torus at each double point of the projection of $L$ onto $\Sigma$; see the left side of Figure~\ref{harer2}. 
\begin{figure}[ht]
  \relabelbox \small {
  \centerline{\epsfbox{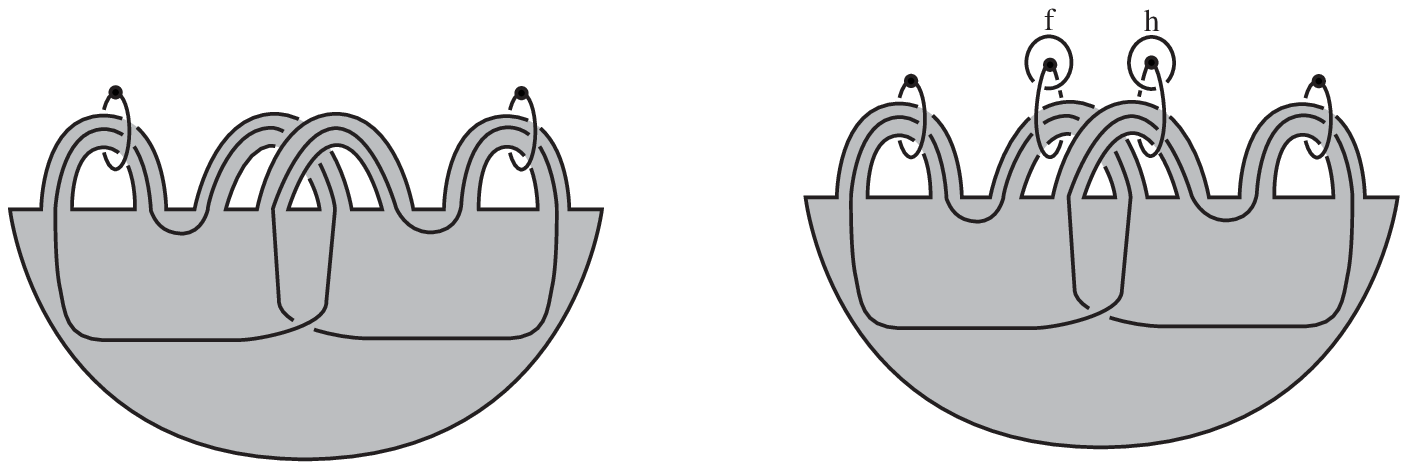}}} 
   \relabel {f}{$\pm 1$}
   \relabel {h}{$\pm 1$}
  \endrelabelbox
        \caption{}
        \label{harer2}
\end{figure}
To arrange that the new $\Sigma$ is a fiber in a non-tivial Lefschetz fibration on $X_1$, we add canceling pairs of 1- and 2-handles (using $\pm 1$-framed 2-handles) for each newly introduced band, as on the right of Figure~\ref{harer2}. Note that the canceling 2-handle may be isotoped to lie on $\Sigma$. 

We are left to correct the framings on the components of $L.$ 
The operation shown in Figure~\ref{harer3} preserves $X$, and alters the framing of a component $K$ of $L$ by $\pm 1$ with respect to the product framings determined by the old and new surfaces $\Sigma$. 
\begin{figure}[ht]
  \relabelbox \small {
  \centerline{\epsfbox{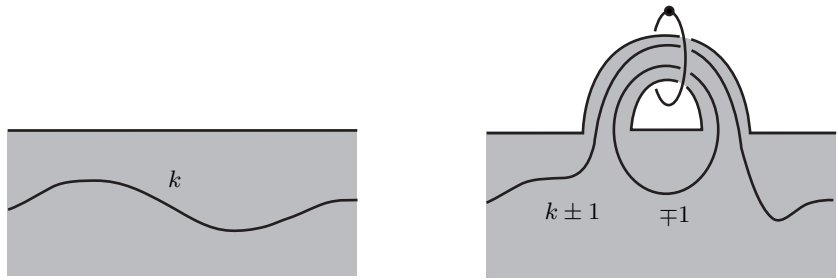}}} 
  \relabel{k}{$k$}
  \relabel {c}{$k\pm 1$}
  \relabel {p}{$\mp 1$}
  \endrelabelbox
        \caption{A stabilization operation for $\Sigma.$ }
        \label{harer3}
\end{figure}
By repeating this procedure as necessary, we can arrange that the framing of $K$ differs from the product framing by either plus or minus one. Moreover, the newly introduced $\mp 1$-framed 2-handles satisfy the framing requirement as well.

We remark that by carefully analyzing this proof one can construct an achiral Lefschetz fibration for $X$ with genus bounded above by the bridge number of the link $L$ (as measured with respect to its projection onto $\Sigma$).
\end{proof}

Let $f\colon\, X\to D^2$ be an achiral Lefschetz fibration with bounded fibers. As above we can describe $X$ as $\Sigma\times D^2$
with 2-handles attached to vanishing cycles $\gamma_1,\ldots
\gamma_k$ with framing one less than the product framing and attached to vanishing cycles $\gamma'_1,\ldots, \gamma'_{k'}$
with framing one more than the product framing. Let $\Sigma'$ be the surface obtained from $\Sigma$ by attaching a 1-handle.
Let $\gamma$ be a simple closed curve embedded in $\Sigma'$ that intersects the cocore of the new 1-handle exactly once. 
A \dfn{positive (negative) stabilization} of this achiral Lefschetz fibration $f$ is the \aLf described as $\Sigma'\times D^2$
with 2-handles attached to $\gamma_1,\ldots,\gamma_k$ and $\gamma'_1,\ldots, \gamma'_{k'}$ as above and a 2-handles attached
to $\gamma$ with framing one less (one more) than the product framing. Note that stabilizing results in an \aLf of the same
4-manifold $X.$ The achiral Lefschetz fibration $f$ induces an open book decomposition $(B,\pi)$ of $\partial X$ and the
positively (negatively) stabilized \aLf also induces an open book decomposition $(B',\pi')$ of $\partial X.$  The open
book $(B',\pi')$ is said to be obtained by positive (negative) stabilization. 

\section{Contact geometry and open book decompositions}\label{sec:obd}
An \dfn{oriented contact structure} on an oriented 3-manifold $M$ is a hyperplane field $\xi$ that
can be written as the kernel of a 1-form $\alpha$ such that $d\alpha$ is non-degenerate when restricted
to $\xi.$ In other words $\xi=\ker \alpha$ and $\alpha\wedge d\alpha\not=0.$ 
We assume the reader is familiar with the basic notions from contact geometry 
(such as Legendrian knot, Thurston-Bennequin invariant and so on).  For a review, see \cite{Etnyre03}.
A contact structure $\xi$ on $M$ is said to be supported by an open book $(B,\pi)$ if $\xi$ is isotopic
to a contact structure given by a 1-form $\alpha$ satisfying $\alpha>0$ on positively oriented tangents
to the binding $B$ and $d\alpha$ is a positive volume form on each page of the open book. Thurston and
Winkelnkemper \cite{ThurstonWinkelnkemper75} showed that any open book 
supports a contact structure. In addition it is fairly simple
to show that two contact structures supported by the same open book are isotopic. Recently Giroux \cite{Giroux02} has
strengthened this connection between contact structures and open book decompositions.
\begin{Thm}[Giroux 2002, \cite{Giroux02}]\label{thm:giroux}
Let $M$ be a closed oriented 3-manifold. There is a one-to-one correspondence between oriented 
contact structures on $M$ up to isotopy and open book decompositions of $M$ up to positive stabilization
(and isotopy).
\end{Thm}

Given a Legendrian knot $L$  let 
$N(L)$ be a {\em standard tubular neighborhood} of the Legendrian curve $L.$ This means the
neighborhood has convex boundary and two parallel dividing curves (see \cite{EtnyreHonda01b}).  Choose a framing for $L$  
so that the meridian has slope $0$ and the dividing curves have slope $\infty$.  With respect to 
this choice of framing, a \dfn{$\pm 1$ contact
surgery} is a $\pm 1$ Dehn surgery, where a copy of $N(L)$ is glued to $M\setminus N(L)$ so that the new 
meridian has slope $\pm 1$.  Even though the boundary characteristic foliations may not 
exactly match up a priori,  we may use Giroux's Flexibility Theorem 
\cite{Giroux91, Honda00a} and the fact that they have the same dividing set to make the 
characteristic foliations agree. This gives us a new contact 
manifold $(M',\xi').$  For a detailed discussion of contact surgery see
\cite{DingGeigesStipsicz04}. The following is a well known theorem, see for example \cite{EtnyreHonda02a}.
\begin{Thm}\label{LegSurgonPage}
Suppose that $L$ is a Legendrian knot in the contact manifold $(M,\xi), \xi$ is supported by the open
book $(B,\phi)$ and $L$ is contained in a page of the open book. The contact manifold
obtained from $(M,\xi)$ by $\pm 1$ contact surgery on $L$ is equivalent to the one compatible with
the open book with monodromy $\phi\circ D^\mp_\alpha.$
\end{Thm}
Returning to Lefschetz fibrations,
let $f\colon\, X\to D^2$ be an achiral Lefschetz fibration with bounded fibers. As at the end of Section~\ref{sec:lf} we can describe $X$ as 
$\Sigma\times D^2$
with 2-handles attached to vanishing cycles $\gamma_1,\ldots
\gamma_k$
and $\gamma'_1,\ldots, \gamma'_{k'}$ with the appropriate framings.
The Lefschetz structure on $\Sigma \times D^2$ induces an open book and 
hence a contact structure on $\partial (\Sigma\times D^2).$ Using the Legendrian realization principle \cite{Honda00a} we
can assume the (non-null homologous) 
vanishing cycles are sitting on the various pages of the open book as Legendrian curves. Moreover the
contact structure induced on $\partial X$ is the one obtained from the contact structure on $\partial (\Sigma\times D^2)$
by $\pm 1$-contact surgeries on the vanishing cycles. 

\begin{figure}[ht]
  \relabelbox \small {\epsfxsize=4.5in
  \centerline{\epsfbox{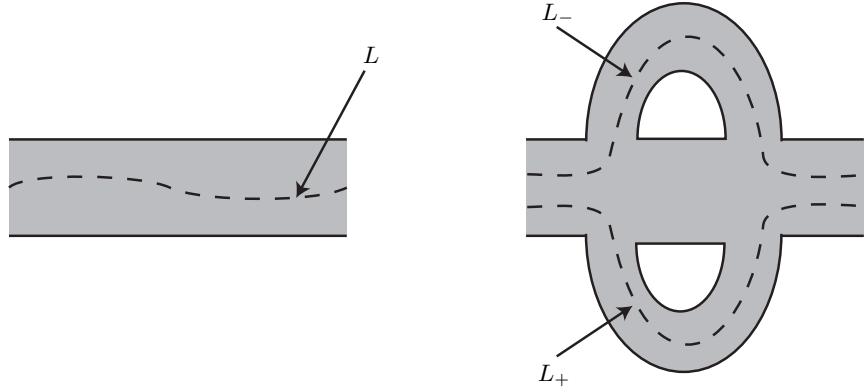}}} 
  \relabel{l}{$L$}
  \relabel {p}{$L_+$}
  \relabel {n}{$L_-$}
  \endrelabelbox
        \caption{A neighborhood of a piece of $L$ in $\Sigma,$ left. ($L$ is oriented so it points towards the left.) 
          The twice stabilized open book, right. If the two stabilizations are positive then $L_\pm=S_\pm(L)$ and
          if the stabilizations are negative then $S_\pm(L_\pm)=L.$}
        \label{stab}
\end{figure}

In our discussion below it will be useful to see how to stabilize (and destabilize) a Legendrian knot on a page of an open
book so that the stabilized knot is also on a page of an open book. Given an oriented Legendrian knot $L$,  let $S_+(L)$ and $S_-(L)$ be the positive and negative stabilizations  of $L$ obtained by adding down or up ``zig-zags''. 
\begin{Lem}\label{lem:stab}
Let $(B,\phi)$ be an open book decomposition supporting the contact structure $\xi$ on $M.$ Suppose
$L$ is a Legendrian knot in $M$ that lies on a page of the open book. If we positively stabilize 
$(B,\phi)$ twice as shown in Figure~\ref{stab} then we may isotop the page of the open book so that
$S_+(L)$ and $S_-(L)$ appear on the page as seen in Figure~\ref{stab}. 

If $(B,\pi)$ is 
negatively stabilized twice as shown in the figure, then the contact structure supported by this
open book is no longer $\xi$ but we still see $L$ as a Legendrian knot in the new contact structure.
Moreover $L_+$ and $L_-$ in the figure are now the positive and, respectively, negative destabilizations of $L.$ That
is $S_\pm(L_\pm)=L.$ 
\end{Lem}
This Lemma is relatively easy to prove, see \cite{EtnyreOBN}.

\section{Overtwisted contact structures and homotopy classes of plain fields.}
Contact structures in dimension three fall into two disjoint classes: tight and overtwisted. A contact
manifold  $(M,\xi)$ is called \dfn{overtwisted} if there is an embedded disk $D$ such that $T_xD=\xi_x$ for all
$x\in \partial D.$ If $\xi$ is not overtwisted it is called tight. 
One may easily prove \cite{EtnyreOBN} that if one negatively stabilizes an open book then the associated contact 
structure is overtwisted.
We have the following fundamental theorem of Eliashberg. 
\begin{Thm}[Eliashberg 1990, \cite{Eliashberg90a}]\label{thm:eliashberg}
If two overtwisted contact structures are homotopic as plane fields then they are isotopic as contact structures. 
\end{Thm}
Using this theorem we can understand overtwisted contact structures by understanding their homotopy class of
plane field.  According to \cite{Gompf98} the homotopy class of an oriented plane field $\xi$ on $M$ is completely 
determined by two invariants. To simplify the discussion we will assume $H^2(M;\Z)$ has no 2-torsion (this
will suffice for our applications).
The first invariant is the first Chern class (a.k.a. Euler class) $c_1(\xi)\in H^2(M;\Z)$, which is simply the obstruction to the existence of a non-zero section of $\xi.$ 
Suppose the contact manifold $(M,\xi)$ is supported by an open book $(B,\pi)$ that is induced as the boundary of the
achiral Lefschetz fibration $f\colon\, X\to D^2.$ We describe this Lefschetz fibration as in Section~\ref{sec:lf} with vanishing cycles
$\gamma_1,\ldots \gamma_k$ and $\gamma'_1,\ldots, \gamma'_{k'}.$  We can assume the $\gamma_i$ and $\gamma'_i$ are Legendrian knots
in $\partial \Sigma\times D^2.$
A slight generalization of a formula from \cite{Gompf98} (see \cite{OzbagciStipsicz04}) computes the Poincar\'e dual to $c_1(\xi)$ as
\begin{equation}\label{c1}
P.D. c_1(\xi)=\sum_{i=1}^k r(\gamma_i) c_i + \sum_{i=1}^{k'} r(\gamma'_i) c'_i,
\end{equation}
where the $c_i$ and $c'_i$ are the image of the cocores of the 2-handles attached to the $\gamma_i$ and $\gamma_i'$'s
under the boundary map in the long exact sequence of the pair $(X,M).$

The second invariant of a homotopy class of oriented plane fields is the
so called 3-dimensional invariant $d_3(\xi)$, which is a rational number well-defined modulo the divisibility of
$c_1(\xi).$  We will only describe how to compute $d_3(\xi)$ when $c_1(\xi)=0.$ To this end let $M$ and 
$X$ be as above.
Then we have 
\begin{equation}\label{d3}
d_3(\xi)= \frac 14 (c^2(X) -3\sigma(X) - 2\chi(X)) + q,
\end{equation}
where $\sigma$ is the signature of $X,$ $\chi$ is the Euler characteristic, and $q$ is the number of negative vanishing 
cycles of $X.$ 
The number $c^2(X)$ is the square of the class $c(X)$ with Poincar\'e dual
\[
\sum_{i=1}^k r(\gamma_i) C_i + \sum_{i=1}^{k'} r(\gamma'_i) C'_i,
\]
where the $C_i$ and $C_i'$'s are the cocores of the 2-handles attached along $\gamma_i$ and $\gamma_i'.$
Note that $c(X)|_M=c_1(\xi)$, which we are assuming to be zero. Thus $c(X)$, 
which naturally lives in $H^2(X;\Z)$, comes from a class in $H^2(X,\partial X;\Z)$ and thus can be squared.
Formula~\eqref{d3} is a slight generalization of the one given in \cite{DingGeigesStipsicz04}, where it was assumed that
$X$ had no 1-handles. Their proof caries over to our case. In particular, according to \cite{Gompf98}, 
$d_3(\xi)=\frac 14 (c^2(Y) -3\sigma(Y) - 2\chi(Y))$ where $Y$ is any almost complex 4-manifold with  
$M=\partial Y$ and $\xi$ is the set of almost complex tangencies to $M.$ If $X$ is as above then
there is a natural almost complex structure on $Y=X\#_q \C P^2$ (see \cite{DingGeigesStipsicz04}) with $\xi$ the set
of complex tangencies. Moreover, 
\[
c_1(Y)=c(X) + (3, \ldots, 3) \in H^2(X;\Z)\oplus_q H^2(\C P^2;\Z),
\]
$\sigma(Y)=\sigma(X)+q,$ and $\chi(Y)=\chi(X)+q.$ The formula follows.

\section{Proof of Theorem~\ref{main}}
We are now ready to establish the existence of Lefschetz fibrations.
\begin{proof}[Proof of Theorem~\ref{main}]
We begin by giving $X$ an arbitrary handlebody structure, letting $X_1$ denote the union of the
$0$-, $1$- and $2$-handles of $X$, and $X_2$ denote the union of the $3$- and $4$-handles. Then as
each $X_i$ is a 2-handlebody, we can use Theorem~\ref{thm:harer} to 
find achiral Lefschetz fibrations (with bounded fibers)
$f_1:X_1\to D^2$ and $f_2:X_2\to D^2$, with each inducing an open book structure on the common
boundary $\bdy X_1=-\bdy X_2$. 
We can stabilize the achiral Lefschetz fibrations so that each 
has fibers with connected boundary.

If these induced open books are the same (under the identification of $\partial X_1$ and $-\partial X_2$ used to reconstruct
$X$), we can attempt to reconstruct $X$ from the pieces $X_1$
and $X_2$ by gluing them along their boundaries in a two step process. We first glue along the pages of
the open books, by forming 
\[
W=X_1\bigcup_{f_1^{-1}(\bdy D^2)=f_2^{-1}(\bdy D^2)} X_2.
\]
We then
have an \aLf with bounded fibers
\[
f_1\cup f_2: W\to S^2.
\] 
Since the fibers in the two achiral Lefschetz fibrations have connected boundary we see 
$\partial W=S^1\times S^2=S^1\times D^2_1\cup S^1\times D^2_2,$ where $S^1\times D_i^2$ is 
$\partial X_i\setminus f_i^{-1}(\partial D^2).$
Gluing $S^1\times D^2_1$ to $S^1\times D^2_2$ in $\partial W$
will yield $X.$ Gluing an $S^1\times D^3$ to $W$ will produce the same result as gluing $S^1\times
D^2_1$ to $S^1\times D^2_2.$ So we see that
\[
X= W\cup S^1\times D^3
\] 
or said another way there is an embedded curve $\gamma$ in $W$ such that $W=X\setminus N$ where $N$
is an open tubular neighborhood of $\gamma.$
 
Notice that in a collar of $\bdy W=S^1\times S^2$ we may express the above \aLf as the projection
$I\times S^1\times S^2\to S^2$. If we now glue in $D^2\times S^2$ so that each $\bdy D^2\times
\{\mathrm{pt.}\}$ matches to $S^1 \times \{\mathrm{pt.}\}$, then the resulting closed manifold has an
\aLf over $S^2.$ Moreover, this manifold is the result of surgering $X$ along the circle $\gamma.$

The theorem is therefore proven once we establish the following proposition.
\end{proof}

\begin{Prop} Let $X$ be a closed, smooth, $4$-manifold. Then we may write $X=Y_1\cup Y_2$ where each $Y_i$
is a 2-handlebody which admits an \aLf over $D^2$ with bounded fibers of the same genus, and with
the induced open books on $\bdy Y_1=-\bdy Y_2$ coinciding.
\end{Prop}

\begin{proof}
Fix a handle decomposition of $X$ and let $Y_1$ be the union of the 0-, 1- and 2-handles and let
$Y_2$ be the union of the 3- and 4-handles. By Harer's Theorem~\ref{thm:harer} we know there are achiral 
Lefschetz fibrations $f_i\colon Y_i\to D^2, i=1,2,$ with bounded fibers. By adding a canceling 
2-handle/3-handle pair to $X$ if necessary we may assume that $-Y_2$ has one 0-handle and an even number $2k$ of 1-handles.
We may then write $-Y_2$ as $\Sigma\times D^2$ 
where $\Sigma$ is a genus $k$ surface
with one boundary component. Written as such, $Y_2$ has an obvious Lefschetz fibration with no singular fibers.
We take $f_2$ to be this fibration. 

Let $\xi_i$ be the contact structure supported by the open book associated to the \aLf $f_i, i=1,2.$ By
Giroux's Theorem~\ref{thm:giroux} these open books associated to the achiral Lefschetz fibrations will be isotopic,
after possible positive stabilization, if the supported contact structures are isotopic. We 
will show how to choose the achiral Lefschetz fibrations so that the associated contact structures $\xi_1$
and $\xi_2$ are isotopic. We begin by showing they are homotopic as plane fields. To this end notice that
$-Y_2$ supports a Stein structure and hence $\xi_2$ is tight. In addition using Equation~\eqref{c1} we see
$c_1(\xi_2)=0.$ 

Let $Y_1'$ denote the union of the 0- and 1-handles of $Y_1.$ Let $K_1,\ldots, K_l\subset \partial Y_1'$ 
be the attaching circles for the 2-handles in $Y_1.$ We know there is a Lefschetz fibration of
$Y_1'$ so that the $K_i$ are Legendrian in the contact structure supported by the induced open book. 
Moreover we can assume the $K_i$ lie on distinct pages of the open book and the attaching framing is $\pm 1$
the page framing (= contact framing). Attaching the 2-handles now gives a natural achiral Lefschetz fibration 
to $Y_1.$ By Equation~\eqref{c1} the Poincar\'e dual of the first Chern class of the induced contact structure $\xi_1$ is
\[
P.D. c_1(\xi_1)=\sum_{i=1}^l r(K_i) c_i,
\]
where $c_i\in H_1(\partial Y_1;\Z)$ is the image of the cocore of the $i^\text{th}$ 2-handle under the boundary
map in the long exact sequence of the pair $(Y_1,\partial Y_1).$ 

In \cite{Gompf98} it was shown that the parity of
$r(K_i)$ is fixed by the contact framing of $K_i$ and the number of 1-handles $K_i$ runs over. 
We claim that it is possible to alter the Lefschetz fibration (and hence the contact structure) so that any integer
with the right parity can be realized as the rotation number of $K_i.$ 
We begin by stabilizing the Lefschetz fibration of $Y_1'$
one time positively and one time negatively. This does not effect the Chern class of the contact structure induced
on $\partial Y_1'$, though the contact structure is different. Thinking of $K_i$ as a Legendrian knot on a
page of this new open book for $\partial Y_1'$, let $K'_i$ be the result of pushing $K_i$ 
over the two new 1-handles in the page. We can assume $K'_i$ is Legendrian. Using Lemma~\ref{lem:stab} we see that the
framings on $K_i$ and $K_i'$ coming from the page are the same. Moreover, using Lemma~\ref{lem:stab}, we can 
choose the stabilizing 1-handles so that  
$r(K_i')=r(K_i)\pm 2.$ Thus by a sequence of stabilizations we can alter the
rotation number of any $K_i$ by any even number. 
Since $P.D. c_1(\xi_1)|_2 = 0$
(see \cite{GompfStipsicz99}),
it now follows that there exists a sequence of alterations of rotation numbers which gives $c_1(\xi_1)=0$.
(To see this, let $G$ denote the subgroup of $H_1(\partial Y_1)$ generated by $c_1,\ldots,c_l$, and note that the subgroup of even elements of $G$ is generated by $2c_1,\ldots,2c_l$. Hence we may write 
$P.D. c_1(\xi_1)=\sum_{i=1}^l a_i (2c_i)$, which combined with Equation~\eqref{c1} gives
$\sum_{i=1}^l (r(K_i)-2a_i) c_i=0$.)

At this point we may assume that $c_1(\xi_1)=c_1(\xi_2)=0.$
The homotopy class of a plane field with $c_1=0$ is determined by the invariant $d_3.$ 
If $\xi$ is supported by an open book and $\xi'$ by the open book obtained by negatively stabilizing one time, then
Equation~\eqref{d3} yields $d_3(\xi')=d_3(\xi)+1.$ By negatively stabilizing the achiral Lefschetz fibration on
$Y_1$ or $Y_2$ we may assume that $d_3(\xi_1)=d_3(\xi_2)$ and thus $\xi_1$ and $\xi_2$ are homotopic as plane fields.
If we now negatively stabilize the achiral Lefschetz fibrations on each of $Y_1$ and $Y_2$ we may assume the
associated contact structures are overtwisted. Eliashberg's Theorem~\ref{thm:eliashberg} allows us to conclude that $\xi_1$
and $\xi_2$ are isotopic contact structures.
\end{proof}

\begin{proof}[Proof of Corollary~\ref{cor1}]
If $X$ be a simply connected 4--manifold and $\gamma$ the curve identified in the proof of Theorem~\ref{main} 
on which surgery produces a manifold $X'$ with an achiral Lefschetz fibration over $S^2.$ It is well-known (see for example
\cite{GompfStipsicz99}) that $X'$ is either $X\#S^2\times S^2$ or $X\# S^2\widetilde{\times}S^2$, with the outcome determined by the framing of $\gamma.$ (The set of framings can be identified with $\pi_1(\mathrm{SO}(3))=\Z_2.$ )
Moreover, since $\gamma\subset M^3\subset X^4$, a framing of $\gamma$ in $M$ gives a framing of
$\gamma$ in $X.$ Recall that $\gamma$ is the binding of an open book and the framing of $\gamma$ in $M$ 
comes from a page $\Sigma$ in the open book. If we positively stabilize the open book twice as shown in
Figure~\ref{framing} 
\begin{figure}[ht]
  \relabelbox \small {\epsfxsize=3.5in
  \centerline{\epsfbox{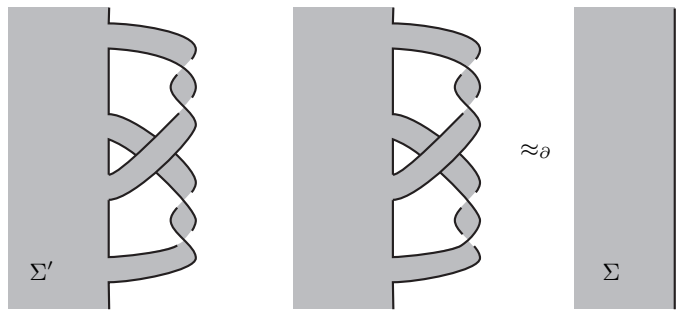}}} 
  \relabel{S}{$\Sigma$}
  \relabel {P}{$\Sigma'$}
  \relabel {e}{$\approx_\partial$}
  \endrelabelbox
        \caption{Surface $\Sigma'$ on the left. In the middle we have $\partial \Sigma'$ with a crossing change and
        	on the right is $\Sigma.$ The boundary of the middle surface is isotopic to the boundary of the right hand
	surface. Although the surfaces are not isotopic, they give the same framing to the knot.}
        \label{framing}
\end{figure}
we have a new knot $\gamma'$ and page $\Sigma'.$ The knots $\gamma'$ and $\gamma$ are homotopic
in $M$ (just change one crossing of $\gamma'$) and hence isotopic in $X.$ The homotopy from $\gamma'$ to $\gamma$
takes the framing on $\gamma'$ coming from $\Sigma'$ to one less than the framing on $\gamma$ coming from $\Sigma.$
Thus the framing on $\gamma'$ in $X$ differs from the framing on $\gamma$ in $X$ and 
therefore surgery on one these curves will yield $X\# S^2\widetilde{\times} S^2$ while surgery on the other will yield $S^2\times S^2$.
\end{proof}

\section{Symplectic and near-symplectic structures}
Let $X$ be a 4-manifold. If we fix a metric $g$ on $X$ we can consider the bundle $\Lambda^2_+$ 
of self-dual 2-forms on $X.$ A closed 2-form $\omega$ on $X$ is 
a \dfn{near-symplectic structure} if  $\omega^2\geq 0$ and there is a metric $g$ such that $\omega$ is harmonic and
transverse to the 0-section of $\Lambda^2_+$
By transversality one can see that the zeros $Z$ of
$\omega^2$ form a union of embedded circles.  Honda \cite{Honda04} showed that each component of $Z$ has 
a neighborhood $S^1\times B^3$ where $\omega$ can be written as one of two models. The ``orientable''
model is $dt\wedge dh + \star_3 dh,$ where
$h$ is a Morse function on $B^3$ with one index 1 or index 2 critical point at 0 and $\star_3$ is the 3-dimensional 
Hodge star operator. The ``non-orientable''  model is a $\Z_2$ quotient of the above model. One may define a near-symplectic
structure without regard to a metric by demanding that $\omega$ is closed and symplectic away from a union of circles $Z$, and
near each component of $Z$ has a model as above. There has been great interest in near-symplectic structures following
work of Taubes \cite{Taubes98} that suggests they might be used to give a ``geometric'' understanding of Seiberg-Witten theory.

If $X$ is allowed to have boundary it is possible that the near-symplectic form $\omega$ degenerates along properly embedded
arcs in $X.$ In this article we will assume that $Z$ is always a union of circles in the interior of $X.$ Given this we can
discuss the convex boundary of a near-symplectic 4-manifold $(X,\omega).$  We say $\partial X$ is \dfn{convex}, or 
strongly convex, 
if there is a vector field $v$ defined near $\partial X,$ transverse to $\partial X$, whose flow expands $\omega,$ namely
\[L_v\omega=c\omega,\]
where $L$ denotes Lie derivative and $c$ is a positive constant.  
The 1-form $\alpha=(\iota_v\omega)_{\partial X}$ is a contact form 
on $\partial X.$ Setting $\xi=\ker \alpha,$ we will say $(X,\omega)$ is a near-symplectic filling of $(\partial X, \xi).$ 
It is a standard fact, see \cite{EtnyreOBN}, that if $\beta$ is any other contact form for $\xi$ then there is a neighborhood of $\partial X$ in 
$X$ that is symplectomorphic to a (one-sided) neighborhood of the graph of some function in $(\partial X)\times \R$ with 
symplectic form $d(e^t \beta),$ where $t$ is the coordinate on $\R.$

The following result, restated for our context, is well-known.
\begin{Thm}[Eliashberg \cite{Eliashberg90a}
and Weinstein \cite{Weinstein91}]\label{build_minus_one}
Suppose $(X,\omega)$ is a near-symplectic filling of $(M,\xi)$, and $L$ is a Legendrian knot in $(M,\xi)$. If a 2-handle is attached along $L$ with contact framing $-1$, then $\omega$ may be extended over the 2-handle to obtain a near-symplectic filling of $(M,\xi')$, where $\xi'$ is the contact structure obtained by contact $-1$ surgery on $L$. Moreover, $(M,\xi')$ is strongly convex. There are no new circles of degeneration in the extended near-symplectic structure.
\end{Thm}

We now turn to establishing a version of this Theorem for $+1$-framed surgeries. To this end we first observe an alternate
description of $+1$-contact surgery. 
\begin{Thm}\label{altsurg}
Let $L$ be a Legendrian knot in a contact manifold $(M,\xi).$ The contact structure obtained from $\xi$ by performing a
$+1$-contact surgery on $L$ is the same as the contact structure obtained by performing a Lutz twist on the positive transverse
push-off of $L$ followed by a Legendrian surgery.
\end{Thm}
We recall the definition of a Lutz twist. If $\gamma$ is a knot transverse to the contact planes of $(M,\xi)$
then $\gamma$ has a standard neighborhood contactomorphic to a neighborhood of the image of the $z$-axis in $\R^2/\sim,$
where $(r,\theta,z)\sim(r,\theta,z+1),$ with contact structure $dz+r^2d\theta.$ We can identify slopes $s\in \R\cup \{\infty\}$
of a linear foliation of $T^2$ by angles $\theta_s\in \R/\pi\Z.$ To distinguish different amounts of twisting we will lift $\theta_s$ to $\R.$
We can express a neighborhood of $\gamma$ as $N=S^1\times D^2$ and assume that on concentric tori $T_a=\{r=a\}$ the characteristic foliation is linear with
monotonically decreasing (as $a$ increases) slope ranging in $(0, -\epsilon].$ (This description uniquely determines the contact structure
on $N.$) If we leave the contact structure $\xi$ the same outside
$N$ but change it so that the slopes of the characteristic foliation on the $T_a$ range in $(0,-\pi-\epsilon]$ then we get a well defined
contact structure $\xi'$ on $M.$ This contact structure is said to be obtained from $\xi$ by a (half-)Lutz twist along $\gamma.$

\begin{proof}
Consider a Legendrian knot $L$ in $(M,\xi).$ Let $N(L)$ be a standard neighborhood of $L.$  (See Section~\ref{sec:obd}.)
We pick a framing on $N(L)$ so that the dividing curves on $\partial N(L)$ have slope $\infty$ and the meridian has slope $0.$
Let $L_+$ be the positive transverse push-off of $L$ contained in $N(L),$  and let $\xi'$ denote the contact structure that is the result of performing a Lutz twist on $L_+.$
The contact structure $\xi'$ agrees with $\xi$ outside $N(L)$ and on $N(L)$ the contact structure $\xi$ is determined
by the fact that the slopes of the characteristic foliation on concentric tori range in $(0,-\frac{\pi}{2}]$, and $\xi'$ by the fact that the slopes
range in $(0,-\frac{3\pi}{2}].$ 
Inside the neighborhood $N(L)$ we may find a standard neighborhood of a Legendrian knot with twisting $2$ in our chosen framing 
as follows.
Break $N(L)$ into two pieces $N_1\cup N_2$ where $N_1$ is a solid tori containing $L_+$ with
slope ranging in $(0,-\frac{5\pi}{6}]$ and $N_2=T^2\times[0,1]$ with slopes ranging in $[-\frac{5\pi}{6}, -\frac{3\pi}{2}].$ (So $N(L)$
is split along a torus with dividing slope $\frac 12.$) The solid torus $N_1$ is a standard neighborhood of a Legendrian knot $L'$ with
twisting number $2$ with respect to the framing chosen on $N(L)$ (see \cite{EtnyreHonda01b}). Thus if we perform a Legendrian surgery on $L'$ this will,
topologically, correspond to a $+1$ surgery on $L.$ Moreover, one can check that the contact structure on $N(L),$ after Legendrian
surgery on $L', $ has slopes of the characteristic foliation on concentric tori ranging in $(-\frac{3\pi}{4}, -\frac{3\pi}{2}].$ Such a contact structure
on a solid torus is tight. Thus we have removed $N(L)$ from $(M,\xi)$ and reglued it with a $+1$-twists and extended $\xi|_{M\setminus N(L)}$
to the surgered manifold so that it is tight on the surgery torus. This is precisely a $+1$-contact surgery.
\end{proof}

To perform a Lutz twist via a near-symplectic cobordism we recall the following result.
\begin{Thm}[Gay and Kirby 2004, \cite{GayKirby04}]
Let $(M,\xi)$ be a contact manifold and
let $\xi'$ be obtained from the contact structure $\xi$ by  a Lutz twist along the transverse curve $\gamma.$ Assume the
Lutz twist occurred in the neighborhood $N$ of $\gamma.$ If $(X,\omega)$ is a near-symplectic filling of $(M,\xi),$
then $\omega$ may be extended over $X\cup M\times[0,1],$ where $\partial X$ and $M\times\{0\}$ are identified, to be
a near-symplectic filling of $(M,\xi').$ Moreover, $\omega$ is symplectic on $(M\setminus N)\times [0,1]$ and $\omega$ 
has one singular circle in $N\times[0,1].$
\end{Thm}
This theorem, coupled with Theorem~\ref{build_minus_one} and 
the proof of Theorem~\ref{altsurg} proves the following result.
\begin{Thm}\label{build_plus_one}
Suppose $(X,\omega)$ is a near-symplectic filling of $(M,\xi)$, and $L$ is a Legendrian knot in $(M,\xi)$. If a 2-handle is attached along $L$ with contact framing $+1$, then $\omega$ may be extended over the 2-handle to obtain a near-symplectic filling of $(M,\xi')$, where $\xi'$ is the contact structure obtained by contact $+1$ surgery on $L$. Moreover, $(M,\xi')$ is strongly convex. There is exactly one singular circle of $\omega$ in the near-symplectic structure on the 2-handle.
\end{Thm}

We now rephrase Theorem 1.1 of \cite{Eliashberg04} for our current purposes.
\begin{Thm}[Eliashberg 2004, \cite{Eliashberg04}]\label{capfiber}
Let $(X',\omega)$ be a near-symplectic filling of $(M,\xi)$ and $(B,\pi)$ be an open book decomposition supporting the
contact structure $\xi.$ Let $X$ be $X'$ with 2-handles attached to $B$ with framing given by the pages of the open book. 
One may extend the near-symplectic structure on $X'$ to $X$ so that no new circles of degeneration are added and
so that the near-symplectic structure restricted to each fiber in  the surface bundle $\partial X$ is symplectic. 
\end{Thm}
The proof in \cite{Eliashberg04} goes through unchanged--in fact, the proof can be simplified since we are assuming $\partial X$ is
strongly convex.

This theorem points out the need to study symplectic surface bundles. 
A \dfn{symplectic bundle over $S^1$}  is a 
3-manifold $M$ that fibers over the
circle together with a closed 2-form $\omega$ which is positive on each fiber. The kernel of
$\omega$ defines a line field that is transverse to the fibers of the fibration. An orientation
on $M$ and on the fibers induces an orientation on the line field and thus we can fix a fiber $\Sigma_0$ of the
fibration and use the line field to define a return map $H_{(M,\omega)}\colon \Sigma_0\to \Sigma_0,$ called the
\dfn{holonomy} of the symplectic fibration. The holonomy $H_{(M,\omega)}$ is a symplectomorphism of $(\Sigma_0,\omega|_{\Sigma_0}).$
If we normalize $\omega$ so that it integrates to 1 on
each fiber of the fibration then the holonomy determines $(M,\omega)$ up to fiber preserving
diffeomorphism. 

Recall a symplectomorphism $f$ of a surface $(\Sigma, \omega)$ is \dfn{Hamiltonian} if there are functions $H_t\colon \Sigma\to \R$
such that the vector field $X_t$ determined by
\[\iota_{X_t}\omega = dH_t\]
generates a flow whose map at $t=1$ is $f.$ A Hamiltonian diffeomorphism is always isotopic to the identity, although the converse is false.
The main result we need is the following.
\begin{Lem}[Eliashberg 2004, \cite{Eliashberg04}]\label{confill}
Suppose the holonomy of the symplectic fibration $(M,\omega)$ is a Hamiltonian diffeomorphism.
Then there is a symplectic form $\Omega$ on $X=\Sigma\times D^2,$ such that $\partial
(X,\Omega)=(M,\omega),$ where $\Sigma$ is the fiber of the fibration. We can moreover assume that
$\Sigma\times\{pt\}$ and $\{pt\}\times D^2$ are symplectic submanifolds in $X.$
\end{Lem}
Since not all symplectomorphisms isotopic to the identity 
are Hamiltonian, we will need a criterion below to determine when a symplectomorphism is Hamiltonian. For this, we consider the flux map. Let $\phi_t, 0\leq t\leq 1,$ be
a path of symplectomorphisms of a surface $(\Sigma, \omega).$ There is a family of vector fields $X_t$ determined by
\[X_t(\phi_t(x))= \frac{d\phi_t(x)}{dt}.\]
The flux of the path $\phi_t$ is defined to be
\[\text{Flux}(\{\phi_t\})=\int_0^1 [\iota_{X_t} \omega]\,  dt,\]
where $[\cdot]$ denotes cohomology class.
This is an element in $H^1(\Sigma;\R).$ (If $\Sigma$ is a torus then the flux is in $H^1(\Sigma;\R/\Z).$)
One can show that the flux only depends on the path of symplectomorphisms up to homotopy with fixed endpoints. 
A path of symplectomorphisms is homotopic to a Hamiltonian path if and only if its flux is zero, see \cite{Calabi70}. There is an
alternate interpretation of flux (see \cite{McDuffSalamon98}) that will be more useful below. Let $(M,\omega)$ be a symplectic bundle over $S^1$
with a holonomy map isotopic to the identity; assume also that the fiber $\Sigma$ has genus greater than 1.
In this case we can identify $M$ as $\Sigma \times S^1$ (up to fiber preserving isotopy). Under this identification
we have a map from $H_1(\Sigma)$ to $H_2(M)$ that sends $[c]$ to $[c\times S^1].$ With this understood the flux
of the holonomy is the map $H_1(\Sigma)\to \R$ given by
\[\text{Flux}([c])= \int_{c\times S^1} \omega.\]

\section{Near-symplectic structures and achiral Lefschetz fibrations}\label{nss}

We are now ready to prove our main theorem concerning near-symplectic structures.

\begin{proof}[Proof of Theorem~\ref{nssonlf}]
Suppose $f\colon X\to S^2$ is an achiral Lefschetz fibration with fiber $\Sigma$ and section $S.$ Let $N_\Sigma$ be
a neighborhood of $\Sigma$ that is fibered by non-singular fibers of $f$ and let $N_S$ be a neighborhood of $S$ that
contains no singular points of $f.$ Let $X'=X\setminus(N_\Sigma\cup N_S).$ We can describe $X'$ as an achiral Lefschetz
fibration over $D^2$ with non-singular fibers $\Sigma'$ by restricting $f$ to $X'.$ Here $\Sigma'$ is simple $\Sigma$ with an
open disk removed.  Thus $X'$ may be built from $\Sigma'\times D^2$ by attaching 2-handles along curves
$\gamma_1,\ldots , \gamma_k$ and $\gamma_1',\ldots, \gamma'_{k'}$ on fibers 
with framing one less and, respectively, one more than the fiber framing.  
We know $\Sigma'\times D^2$ has a symplectic structure with convex boundary, see \cite{EtnyreOBN}. Moreover, the open book
induced on the boundary from the product structure supports the induced contact structure on the boundary. 
Thus from Theorems~\ref{build_minus_one} and \ref{build_plus_one} we see that $X'$ has
a near-symplectic structure, with one degenerate circle for each $\gamma'_i.$ In addition, with this near-symplectic structure
the boundary of $X'$ has a contact structure supported by the open book induced by $f.$

Notice that $N_\Sigma\cup N_S$ is simply $\Sigma\times D^2$ union a 2-handle $h.$ If we view $h$ as attached to $X'$ instead of $N_\Sigma,$
it will be attached to the binding of the open book with framing coming from the fibers of the open book. This process simply caps off the fibers $\Sigma'$ to recover the fibers $\Sigma,$ since the result of moving $h$ to $X'$ is $X\setminus (\Sigma\times D^2)$
which has an achiral Lefschetz fibration over $D^2$ with fiber $\Sigma.$ Thus the result of surgery along the binding of 
the open book for $\partial X'$ with framing coming from the fiber is $\Sigma\times S^1.$ Moreover, by Lemma~\ref{capfiber}, we know the near-symplectic form on $X'$ extends to $X'\cup h$ so that the surfaces $\Sigma\times\{pt\}$ are all symplectic.
If the holonomy of this symplectic fibration is trivial, or Hamiltonian isotopic to the identity,  
then, by Lemma~\ref{confill},  we may extend this near-symplectic structure over $N_\Sigma$ thus 
constructing the near-symplectic structure on $X.$ We are left to show that we can arrange the holonomy to be trivial.

We begin by observing that the monodromy of the open book is a composition of Dehn twists parallel to the boundary of
$\Sigma'.$ (Of course the monodromy expressed in terms of Dehn twists along the $\gamma_i$'s and $\gamma'_i$'s
might look more complicated, but it will be isotopic to this.) 
We may assume that the monodromy is supported in a collar neighborhood of the boundary, and we write $\Sigma''$ for the complement of this neighborhood in $\Sigma'$. The complement in $M=\partial X'$ of the binding and the support of the monodromy (in all fibers) can be written $\Sigma''\times S^1,$ and is denoted $M'.$ 
The contact structure $\xi$ on $M'$ is isotopic to one
given by the kernel of $\alpha|_{M'} = K\, dt + \lambda,$ where $K$ is any large positive constant, $t$ is the coordinate on $S^1$
and $\lambda$ a primitive for a volume form on $\Sigma''.$ It is easy to see the Reeb vector field for $\alpha$ is 
$X=\frac{\partial}{\partial t}.$ Now if we consider the 4-manifold $Y=M\times[a,b]$ with symplectic form $\omega=d(e^s \alpha),$
where $s$ is the coordinate on the interval factor, the upper boundary of $Y$ is convex and induces the contact structure $\xi.$
In addition, the kernel of $\omega|_{M\times\{a\}}$ is spanned by the Reeb vector field and the flow of the Reeb vector field
induces the identity return map on the $\Sigma''$ part of a page of the open book. Now if we attach a 2-handle to 
$Y$ along the binding of the open book in $M\times\{b\}$ as in Lemma~\ref{capfiber}, then we obtain a symplectic manifold $Y'$ with an upper
boundary $\Sigma\times S^1$, which has symplectic fibers $\Sigma\times\{pt\}.$ Since the symplectic structure is only affected near the attaching region 
for 2-handle, the kernel of the symplectic form restricted to the upper boundary will still induce the identity  map on the
$\Sigma''$ part of the fiber. Given any primitive homology class $h\in H_1(\Sigma;\Z)$ we can represent it by an embedded curve
$c$ contained in $\Sigma''.$ Now $\omega$ restricted to $c\times S^1$ in the upper boundary of $Y'$ is zero (since the 
$\{pt\}\times S^1$ is in the kernel of $\omega|_{\partial Y'}$). Thus $\int_{c\times S^1} \omega =0.$ Using our second interpretation
of flux we see the flux of the holonomy is zero and hence the holonomy map is Hamiltonian isotopic to the identity.

Since $\partial X'$ is strongly convex (and the near-symplectic structure is symplectic there) a neighborhood of $\partial X'$ is 
symplectomorphic to a neighborhood of the graph of a function $g$ in $M\times\R$ with symplectic structure $\omega=d(e^s \alpha).$ 
Let $b$ be any number larger than the maximum of $g$. We can add the collar $\{(x,s)\in M\times \R : g(x)\leq s\leq b\}$ to $X'$ 
and extend the symplectic structure over it so that a neighborhood of $\partial X'$ is symplectomorphic to $M\times[a,b]$, as in the previous
paragraph. We may now attach the 2-handle $h$ to $X'$ and see the holonomy is Hamiltonian isotopic to the identity, as described above.

Finally we note that the section of $X'$ union $h$ given by the co-core of the 2-handle is symplectic. By Lemma~\ref{confill}, this 
section may be symplectically extended over $N_\Sigma$, showing that the original section of $X$ over $S^2$ is symplectic. If we started with more than one section
of $X$ then we could have removed $N_\Sigma$ and neighborhoods of each of these sections to form $X'.$ The argument
above applies equally well to this case.
\end{proof}

\bibliography{references}
\bibliographystyle{plain}
\end{document}